\documentclass[10pt]{sig-alternate}

\usepackage{bm}
\usepackage{epsfig}
\usepackage{graphicx}
\usepackage{subfigure}
\usepackage{float}
\usepackage{cite}
\usepackage{url}
\usepackage{color}
\usepackage{balance}
\usepackage{mdwlist}
\usepackage{multirow}
\usepackage{threeparttable}
\usepackage{enumitem}
\usepackage{amsmath}
\usepackage{stmaryrd}
\usepackage{booktabs}
\usepackage{siunitx}
\usepackage[ruled,vlined]{algorithm2e}
\usepackage{threeparttable}
\usepackage{graphicx}
\usepackage{epsfig,amsmath,amsfonts,cite}
\DeclareMathAlphabet\mathbfcal{OMS}{cmsy}{b}{n}
\newcommand{\ten}[1]{\mathbfcal{#1}}
\newcommand{\mat}[1]{\mathbf{#1}}

\newcommand{\parm}{{\xi}}
\newcommand{\vecpar}{\boldsymbol{\parm}}
\newcommand{\veceta}{\boldsymbol{\eta}}

\newcommand{\parNum}{d}

\newcommand{\out}{y}

\newcommand{\Phimat}{\boldsymbol{\Phi}}
\newcommand{\multiGPC}{\Psi }

\newcommand{\polyInd}{\alpha}
\newcommand{\basisInd}{\boldsymbol{\polyInd}}
\newcommand{\pcOrder}{p}

\newcommand{\yPC}{\sum\limits_{|\basisInd|=0}^{\pcOrder} {c_{\basisInd}  \multiGPC_{\basisInd}  (\vecpar)} }

\newcommand{\muvec}{\boldsymbol{\mu} }
\newcommand{\Sigmamat}{\mathbf{\Sigma}}

\newcommand{\bvec}{\mathbf{b} }

\newcommand{\yvec}{\mathbf{y} }
\newcommand{\cvec}{\mathbf{c} }

\newcommand{\xvec}{\mathbf{x} }
\newcommand{\Amat}{\mathbf{A} }
\newcommand{\Bmat}{\mathbf{B} }
\newcommand{\Cmat}{\mathbf{C} }

\newcommand{\Mmat}{\mathbf{M} }
\newcommand{\Lmat}{\mathbf{L} }
\newcommand{\Gmat}{\mathbf{G} }
\newcommand{\Atensor}{\mathbfcal{A} }

\newcommand{\Btensor}{\mathbfcal{B} }

\newcommand{\reff}[1]{(\ref{#1})}
\def\sssp{\def\baselinestretch{0.88}\large\normalsize}\sssp

 \makeatletter
 \def\@copyrightspace{\relax}
 \makeatother

\begin{document}

\title{Uncertainty Quantification of Electronic and Photonic ICs
with Non-Gaussian Correlated Process Variations}

  \author{
  \alignauthor Chunfeng Cui and Zheng Zhang\\
         \affaddr{Department of Electrical \& Computer Engineering, University of California, Santa Barbara, CA 93106}\\
        \email{{ chunfengcui@ucsb.edu; zhengzhang@ece.ucsb.edu}}
  }

\maketitle

\begin{abstract}
Since the invention of generalized polynomial chaos in 2002, uncertainty quantification has impacted many engineering fields, including variation-aware design automation of integrated circuits and integrated photonics. Due to the fast convergence rate, the generalized polynomial chaos expansion has achieved orders-of-magnitude speedup than Monte Carlo in many applications. However, almost all existing generalized polynomial chaos methods have a strong assumption: the uncertain parameters are mutually independent or Gaussian correlated. This assumption rarely holds in many realistic applications, and it has been a long-standing challenge for both theorists and practitioners.

This paper propose a rigorous and efficient solution to address the challenge of non-Gaussian correlation. We first extend generalized polynomial chaos, and propose a class of smooth basis functions to efficiently handle non-Gaussian correlations. Then, we consider high-dimensional parameters, and develop a scalable tensor method to compute the proposed basis functions. Finally, we develop a sparse solver with adaptive sample selections to solve high-dimensional uncertainty quantification problems. We validate our theory and algorithm by electronic and photonic ICs with 19 to 57 non-Gaussian correlated variation parameters. The results show that our approach outperforms Monte Carlo by $2500\times$ to $3000\times$ in terms of efficiency. Moreover, our method can accurately predict the output density functions with multiple peaks caused by non-Gaussian correlations, which is hard to handle by existing methods.

Based on the results in this paper, many novel uncertainty quantification algorithms can be developed and can be further applied to a broad range of engineering domains.


\end{abstract}

%
%
%
%


\section{Introduction}
Uncertainties are unavoidable in almost all engineering fields, and they should be carefully quantified and managed in order to improve design reliability and robustness. In semiconductor chip design, a major source of uncertainty is the fabrication process variations. Process variations are  significant in deeply scaled electronic integrated circuits (ICs)~\cite{variation2008} and MEMS~\cite{agarwal2009stochastic}, and they have also become a major concern in emerging design technologies such as integrated photonics~\cite{zortman2010silicon}. A popular uncertainty quantification method is Monte Carlo~\cite{MCintro}, which is easy to implement but has a low convergence rate. In recent years, various stochastic spectral methods (e.g., stochastic Galerkin \cite{ghanem1991stochastic}, stochastic testing \cite{zzhang:tcad2013} and stochastic collocation \cite{xiu2005high}) have been developed and have achieved orders-of-magnitude speedup than Monte Carlo in vast applications. These methods represent a stochastic solution as the linear combination of some basis functions (e.g., generalized polynomial chaos~\cite{gPC2002}), and they can obtain highly accurate solutions at a low computational cost when the parameter dimensionality is not high.

Stochastic spectral methods have been successfully applied in the variation-aware modeling and simulation of many devices and circuits, including (but not limited to) VLSI interconnects \cite{Tarek_DAC:08, Wang:2004,Shen2010,cmpt2012}, nonlinear
ICs~\cite{Strunz:2008, Tao:2007, zzhang:tcad2013}, MEMS~\cite{zzhang_cicc2014,agarwal2009stochastic} and photonic circuits \cite{twweng:optsEx}. A major challenge of stochastic spectral methods is the curse of dimensionality: a huge number of basis functions and simulation samples may be required as the number of random parameters becomes large. In recent years, there has been significant improvement to address this challenge. Representative techniques include (but are not limited to) compressive sensing \cite{xli2010, hampton2015compressive}, analysis of variance~\cite{ma2010adaptive,zzhang_cicc2014}, stochastic model order reduction~\cite{el2010variation}, hierarchical methods \cite{zhang2014calculation,zzhang_cicc2014,zhang2015enabling} and tensor computation~\cite{zhang2017big, zhang2015enabling,zhang2017tensor}.

{\bf Major Challenge.} Despite their great success, existing stochastic spectral methods are limited by a long-standing challenge: the generalized polynomial-chaos basis functions require all random parameters to be {\bf mutually independent}~\cite{gPC2002}. This is a very strong assumption, and it fails in many realistic cases. For instance, a lot of device geometric or electrical parameters are highly correlated because they are influenced by the same fabrication steps; many circuit-level performance parameters used in system-level analysis depend on each other due to the network coupling and feedback. Data-processing techniques such as principal or independent component analysis~\cite{wold1987principal,singh2006statistical} can handle Gaussian correlations, but they cause huge errors in general {\bf non-Gaussian correlated} cases. A modified and non-smooth chaos representation was proposed in~\cite{soize2004physical}, and it was applied to the uncertainty analysis of silicon photonics~\cite{twweng:optsEx}. However, the method in~\cite{soize2004physical} does not converge well, and designers cannot easily extract statistical information (e.g., mean value and variance) from the solution.
  \begin{figure*}[t]
	\centering
		\includegraphics[width=5.6in, height=1.8in]{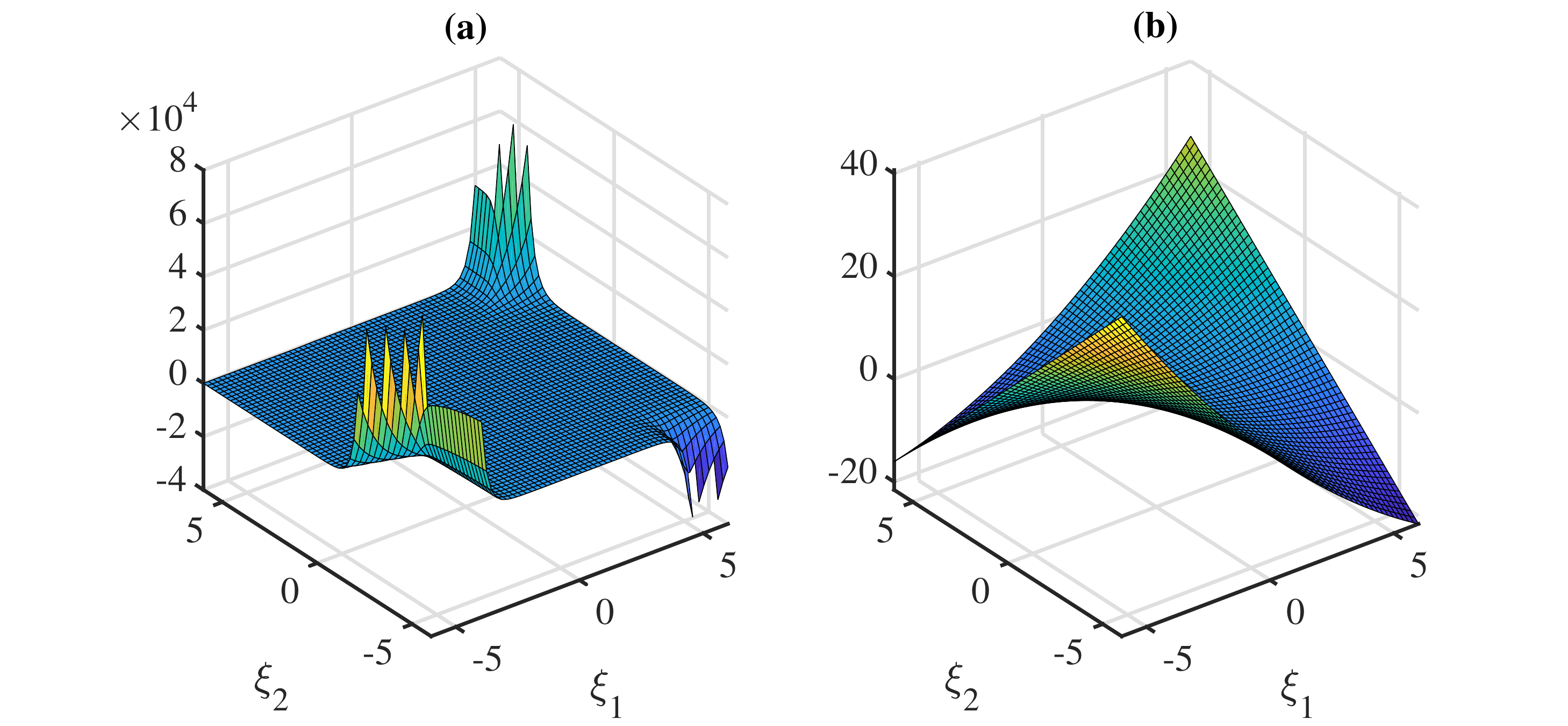}
\caption{ (a): A two-variable basis function by~\cite{soize2004physical}; (b): a basis function obtained by our proposed method.}
	\label{fig:basisfun}
\end{figure*}

{\bf Paper Contributions.} This paper proposes a novel and rigorous solution to handle the challenging non-Gaussian correlated process variations in electronic ICs and integrated photonics. The specific contributions include:


\begin{itemize}[leftmargin=*]
\item Derivation and implementation of a class of basis functions for non-Gaussian correlated random parameters. The proposed basis functions can overcome the theoretical limitations of~\cite{soize2004physical}: they are smooth and can provide highly accurate solutions for non-Gaussian correlated cases (like the standard generalized polynomial chaos~\cite{gPC2002} does for independent parameters), and it can provide closed-form expressions for the  mean values and variance of a stochastic solution. In order to make our methods scalable, we also propose a highly efficient functional tensor-train method to compute the basis functions for many non-Gaussian correlated random parameters equipped with a Gaussian-mixture density function.


\item An adaptive-sampling sparse solver. In order to apply our method to electronic and photonic ICs with  non-Gaussian correlated variations, we develop an $\ell_0$-minimization framework to compute the sparse coefficients of our basis functions. Our {\it main contribution} is an adaptive sampling approach: instead of setting up a compressive-sensing problem using random samples (as done in~\cite{xli2010}), we select the most informative samples via a rank-revealing QR procedure and use a D-optimal method to add new samples and to update the solution.

\item Validation on electronic and photonic ICs. We demonstrate the effectiveness of our framework on electronic and photonic IC examples with $19$ to $57$ non-Gaussian correlated process variations. Our method can accurately predict the statistical information (e.g., multi-peak probability density function and mean value) of the circuit performance, and it is faster than Monte Carlo by about $3000\times$  when the similar level of accuracy is required.
\end{itemize}

Our uncertainty quantification framework has the following two excellent features simultaneously: it does not need any error-prone de-correlation step such as independent component analysis~\cite{singh2006statistical}, and it has the similarly high performance for non-Gaussian correlated uncertainties as generalized polynomial chaos~\cite{gPC2002} does for independent uncertainties.



\section{Preliminaries}
\label{sec:preliminaries}

\subsection{Generalized Polynomial Chaos}
\label{subsec:uq}
Let $\vecpar=[\parm_1, \ldots, \parm_{\parNum}] \in \mathbb{R}^{\parNum}$ denote $d$  random parameters with a joint probability density function $\rho(\vecpar)$, and $\out (\vecpar) \in \mathbb{R}$ be a parameter-dependent performance metric (e.g., the power consumption or frequency of a chip). When $\out (\vecpar)$ is smooth and has a bounded variance, stochastic spectral methods aim to approximate $\out(\vecpar)$ via a truncated generalized polynomial-chaos expansion~\cite{gPC2002}:
\begin{equation}
\label{eq:ygpc}
\out (\vecpar) \approx \yPC,
\end{equation}
where $c_{\basisInd}$ is the coefficient, and $\{{\multiGPC}_{\basisInd} \left(\vecpar\right)\}$ are orthonormal polynomials satisfying
\begin{equation}
 \mathbb{E}\left[{\multiGPC}_{\basisInd} \left( \vecpar \right) \multiGPC_{\boldsymbol{\beta }}\left( \vecpar \right)\right ]=
\left\{\begin{array}{cc}
1,&{\rm\  if\ } \basisInd=\boldsymbol{\beta };\\
0,&\rm{otherwise.}
\end{array}\right.
\end{equation}
Here the operator $\mathbb{E}$ denotes expectation;
$\basisInd=[\alpha_1,\ldots, \alpha_{\parNum}] \in \mathbb{N}^{\parNum}$ is a vector, with each element $\alpha_i$ being the highest polynomial order in terms of $\xi_i$.
The total polynomial order $|\basisInd|=|\alpha_1|+\ldots +|\alpha_d|$ is bounded by $p$, and thus the total number of basis functions is
$N=(p+d)!/(p!d!)$. The unknown coefficients $c_{\basisInd}$'s can be computed via various numerical solvers such as stochastic Galerkin \cite{ghanem1991stochastic}, stochastic testing~\cite{zzhang:tcad2013} and  stochastic collocation~\cite{xiu2005high}. Once $c_{\basisInd}$'s are computed, the mean value, variance and density function of $\out (\vecpar)$ can be easily obtained.

The generalized polynomial-chaos theory~\cite{gPC2002} assumes that all random parameters are mutually independent. In other words, if $\rho _k(\xi _k)$ denotes the marginal density of $\xi_k$, then the joint density is $\rho(\vecpar) = \prod\limits_{k=1}^d {\rho _k(\xi _k)}$. Under this assumption, a multivariate basis function has a product form:
\begin{equation}
\label{eq:gPC_basis}
\multiGPC_{\basisInd}(\vecpar) =\prod\limits_{k=1}^d {\phi_{k,\alpha_k} ( {\xi_k } )}.
\end{equation}
Here ${\phi_{k,\alpha_k} ( {\xi_k } )}$ is a univariate degree-$\alpha_k$ orthonormal polynomial of parameter $\xi_k$, and it is adaptively chosen based on $\rho _k(\xi _k)$ via the three-term recurrence relation~\cite{Walter:1982}.


\subsection{Existing Solutions for Correlated Cases}
The random parameters $\vecpar$ are rarely guaranteed to be independent in realistic cases. It is easy to de-correlate Gaussian correlated random parameters via principal or independent component analysis~\cite{wold1987principal,singh2006statistical}, but de-correlating non-Gaussian correlated parameters can be error-prone. In \cite{soize2004physical}, Soize and Ghanem suggested the following basis function:
\begin{equation}
\label{eq:apc}
\multiGPC_{\basisInd}(\vecpar) =\left(
\frac{\prod\limits_{k=1}^d {\rho _k(\xi _k)}}{\rho(\vecpar)}
\right)^{\frac12}\prod\limits_{k=1}^d {\phi_{k,\alpha_k} ( {\xi_k } )}.
\end{equation}
The modified basis functions are guaranteed to be orthonormal even if
$\rho(\vecpar) \neq \prod\limits_{k=1}^d {\rho _k(\xi _k)}$, but they have two limitations as shown by the numerical results in~\cite{twweng:optsEx}:
\begin{itemize}[leftmargin=*]
\item Firstly, the basis functions are very non-smooth and numerically unstable due to the first part on the right-hand side of \eqref{eq:apc}. This is demonstrated in Fig.~\ref{fig:basisfun}~(a). As a result, the modified basis functions have a much slower convergence rate compared with the standard method in~\cite{gPC2002}.

\item Secondly, the basis functions in \reff{eq:apc} do not allow an explicit expression for the expectation and variance of $\out (\vecpar)$. This is caused by the fact the basis function indexed by ${\basisInd}=0$ is not a constant.
\end{itemize}

\subsection{Background: Tensor Train Decomposition}
\label{subsec:tensortrain}
A tensor is a generalization of a vector and a matrix. A vector $\mat{a} \in \mathbb{R}^n$ is a one-way data array; a matrix $\mat{A} \in \mathbb{R}^{n_1 \times n_2}$ is two-way data array; a tensor $\ten{A} \in \mathbb{R}^{n_1\times n_2 \times \cdots n_d}$ is a $d$-way data array.
We refer readers to~\cite{tensor:suvey,zhang2017tensor} for detailed tensor notations and operations, and its application in EDA~\cite{zhang2017tensor}.

A high-way tensor has $O(n^d)$ elements, leading to a prohibitive computation and storage cost. Fortunately, realistic data can often be factorized using  tensor decomposition techniques~\cite{tensor:suvey}. Tensor train decomposition~\cite{oseledets2011tensor} is very suitable for factorizing high-way tensors, and it only needs $O(dr^2n)$ elements to represent a high-way data array. Specifically, given a $d$-way tensor $\Atensor$, the tensor-train decomposition represents each element $a_{i_1i_2\cdots i_d}$ as
\begin{equation}
a_{i_1 i_2\cdots i_d}=\Amat_1(i_1)\Amat_2(i_d)\ldots \Amat_d(i_d), \forall i_k=1,2,\cdots, n_k,
\end{equation}
where $\Amat_k(i_k)$ is an $r_{k-1}\times r_k$ matrix, and
$r_0=r_d=1$. Given two tensors $\Atensor$ and $\Btensor$ and their tensor train decompositions.
If we want to compute the tensor train decomposition of
their Hadamard (element-wise) product
\begin{equation*}
\mathcal{C}=\Atensor\circ \Btensor \;\; \Longleftrightarrow \;\;c_{i_1 i_2\cdots i_d}=a_{i_1 i_2\cdots i_d}b_{i_1 i_2\cdots i_d},
\end{equation*}
then the result can be directly obtained via
\begin{equation}
\label{equ:TTkron}
\Cmat_k(i_k) = \Amat_k(i_k)\otimes \Bmat_k(i_k).
\end{equation}
Here $\otimes$ denotes a matrix Kronecker product.

\begin{figure*}[t]
	\centering
		\includegraphics[width=6.4in, height=1.7in]{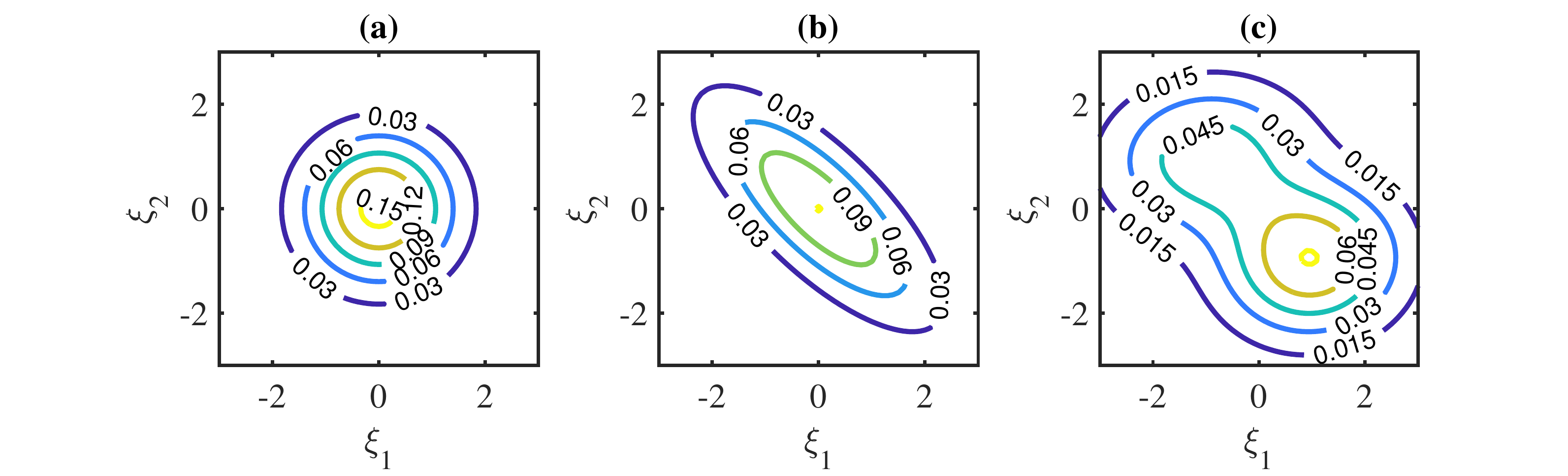}
\caption{Joint density functions for different cases. (a): independent Gaussian; (b): correlated Gaussian; (c): correlated non-Gaussian (e.g., a Gaussian-mixture distribution).}
	\label{fig:gmdistribution}
\end{figure*}

\section{Basis Functions for Non-Gaussian Correlated Cases}
\label{sec:nongaussian}

When process variations are non-Gaussian correlated, the joint probability density of $\vecpar$ is not the product of $\rho_k(\xi_k)$'s, and the basis functions in~\eqref{eq:gPC_basis} cannot be employed. This section derives a set of   multivariate polynomial basis functions. These basis functions can be obtained if a multivariate moment computation framework is available. A broad class of non-Gaussian correlated parameters are described by Gaussian mixture models. For these cases, we propose a fast functional tensor-train method to compute the desired multivariate basis functions.

 \subsection{Proposed Multivariate Basis Functions}
We aim to generate a set of multivariate orthonormal polynomials with respect to the joint density $\rho(\vecpar)$, such that they have the excellent properties of the generalized polynomial chaos~\cite{gPC2002} even for non-Gaussian correlated cases. Several  orthogonal polynomials exist for a few specific density functions~\cite{ismail2017review}. In general, one may construct   multivariate orthogonal
polynomials via the three-term recurrence
in~\cite{xu1993multivariate} or~\cite{barrio2010three}.
However, the theories in ~\cite{xu1993multivariate, barrio2010three} either are hard to implement or can only guarantee week orthogonality.

Inspired by \cite{Golub:1969}, we present a simple yet efficient method for computing  a set of multivariate orthonormal polynomial basis functions.
Let $\vecpar^{\basisInd}=\xi_1^{\alpha_1}\xi_2^{\alpha_2}\ldots\xi_d^{\alpha_d}$ be a monomial indexed by $\basisInd$, then the corresponding moment is
 \begin{equation}
 \label{equ:moments}
 m_{\basisInd} = \mathbb{E}[\vecpar^{\basisInd}]:=\int \vecpar^{\basisInd} \rho(\vecpar) d\vecpar.
 \end{equation}
We intend to construct $N=\frac{(p+d)!}{p!d!}$ multivariate orthonormal polynomials $\{\multiGPC_{\basisInd}(\vecpar)\}$ with their total degrees $|\basisInd|\leq p$. For convenience, we resort all monomials in the graded lexicographic order, denoted as $\bvec(\vecpar)=[b_1(\vecpar),\ldots,b_{N}(\vecpar)]^T$. We further denote the  multivariate moment matrix as $\Mmat$
 \begin{equation}
 \label{equ:mommat}
\mat{M}=\mathbb{E} \left( \mat{b}(\vecpar) \mat{b}^T(\vecpar)\right) \Longleftrightarrow m_{ij}=\mathbb{E}[b_i(\vecpar)b_j(\vecpar)].
 \end{equation}
 For instance, if $d=2$ and $p=2$, then the monomials are ordered as
 \begin{equation*}
 \mat{b}(\vecpar)=[1, \xi_1,\xi_2, \xi_1^2, \xi_1\xi_2, \xi_2^2]^T.
 \end{equation*}
The total number of monomials is $N=6$, and
the corresponding $\mat{M}$ is a $6$-by-$6$ matrix.

Because $\mat{M}$ is a symmetric positive definite matrix, a lower-triangular matrix $\mat{L}$ is calculated via the Cholesky factorization $\Mmat=\Lmat\Lmat^T$. Finally, we define our basis functions as
 \begin{equation}
 \label{equ:basis}
  \boldsymbol{\Psi}(\vecpar) =\Lmat^{-1}\mathbf{b}(\vecpar).
 \end{equation}
Here the $N$-by-1 functional vector $\boldsymbol{\Psi}(\vecpar)$ stores all basis functions $\{\multiGPC_{\basisInd}(\vecpar)\}$ in the graded lexicographic order.

{\bf Properties of the Basis Functions.} Our proposed basis functions have the following excellent properties:
\begin{enumerate}[leftmargin=*]
\item[1).] Our basis functions are smooth (c.f. Fig.~\ref{fig:basisfun} (b)). In fact, the basis function $\multiGPC_{\basisInd}(\vecpar)$  is a multivariate polynomial of a total degree $|\basisInd|$.
It differs from the standard generalized polynomial chaos \cite{gPC2002} in the sense that our basis functions are not the product of univariate polynomials.
\item[2).] All basis functions are orthonormal to each other. This can be easily seen from
\begin{equation}
 \nonumber
  \mathbb{E}\left(\boldsymbol{\Psi}(\vecpar) \boldsymbol{\Psi}^T(\vecpar)\right)=\Lmat^{-1}\mat{M}\mat{L}^{-T}=\mat{I}.
 \end{equation}
This property is important for the sparse approximation and for extracting the statistical information of $\out(\vecpar)$.
\item[3).]  Due to the orthonormality of $\{ \boldsymbol{\Psi}_{\basisInd} (\vecpar)\}$, there exist closed-form formula for the expectation and variance of $\out(\vecpar)$:
\begin{align}
\label{equ:mean}
\mathbb{E}[\out (\vecpar)] &\approx  \sum_{|\basisInd|=0}^p c_{\basisInd} \mathbb{E}\left[\multiGPC_{\basisInd}(\vecpar)\right]=c_{\mat{0}},\\
\text{var}[\out (\vecpar)] &= \mathbb{E} \left(\out^2 (\vecpar) \right)- \mathbb{E}^2 \left(\out (\vecpar) \right) \approx \sum_{|\basisInd|=1}^p c_{\basisInd}^2.
\end{align}

\end{enumerate}

A key step in constructing our basis functions is to compute a set of moments $\{ m_{\basisInd} \}$ for all $|\basisInd|\le 2p$. In general, this can be done with the help of the Rosenblatt transform~\cite{rosenblatt1952remarks}.
Suppose $\rho(\vecpar)=\sum_{i=1}^n w_i\rho_i(\vecpar)$. Based on the joint cumulative density function, the  Rosenblatt formulate a function transformation $\vecpar = \mathcal{T}_i(\veceta)$ such that
$\tilde{\rho}_i(\veceta) =  \beta_i\rho_i( \mathcal{T}_i(\veceta))$ is the density function of the mutually independent parameters $\veceta$.
Here, $\beta_i$ is a coefficient to ensure that the integral of $\tilde{\rho}_i(\veceta)$
is one. With the Rosenblatt transform, we have
\begin{align}
\label{equ:xi2eta_moment}
\nonumber     \mathbb{E}[\vecpar^{\basisInd}] =&\sum_{i=1}^n w_i\int_{-\infty}^{\infty} \vecpar^{\basisInd} \rho_i(\vecpar)  d\vecpar=\sum_{i=1}^n w_i\int_{-\infty}^{\infty} (\mathcal{T}_i(\veceta))^{\basisInd} \tilde{\rho}_i(\veceta) d\veceta.
\end{align}
Then we can use the sparse grid technique~\cite{Gerstner:1998} to numerically compute $\mathbb{E}[\vecpar^{\basisInd}]$. We can also compute a high-dimensional integration via tensor trains as has been done in~\cite{zhang2015enabling}.

 \subsection{Moments for Gaussian-Mixture Models}
\label{subsec:moments}
In practice, semiconductor foundries usually have a lot of measurement data about process variations, and they generate a joint density function $\rho(\vecpar)$ to fit the measurement data set. An excellent choice for this data-driven modeling flow is the Gaussian-mixture model. A Gaussian mixture model   describes the joint density function as
 \begin{equation}
 \rho(\vecpar) =\sum_{i=1}^n w_i \mathcal{N}(\vecpar | \muvec_i, \Sigmamat_i), \; {\rm with}\; w_i>0 , \; \sum_{i=1}^n w_i=1.
 \end{equation}
Here $\mathcal{N}(\vecpar | \muvec_i, \Sigmamat_i)$ is multi-variate Gaussian density function with mean $\muvec_i \in\mathbb{R}^d$ and a positive definite covariance matrix $\Sigmamat_i \in \mathbb{R}^{d \times d}$. Please note that a Gaussian mixture model describes a non-Gaussian correlated joint density function, as shown in Fig.~\ref{fig:gmdistribution}. Now the moment is
\begin{equation}
m_{\basisInd}=\sum\limits_{i=1}^n w_i q_{\basisInd,i}, \; {\rm with}\; q_{\basisInd,i}= \int_{-\infty}^{\infty} \vecpar^{\basisInd} \mathcal{N}(\vecpar | \muvec_i, \Sigmamat_i) d\vecpar,
\end{equation}
and we need to compute $q_{\basisInd,i}$ for $i=1,2,\ldots,n$.

For simplicity, we ignore the index $i$ in $\muvec_i$, $\Sigmamat_i$ and $q_{\basisInd,i}$. Let $\Amat$ be the lower triangular matrix from the Cholesky decomposition of $\Sigmamat$ (i.e.,  $\Sigmamat=\Amat\Amat^T$), and let $\vecpar=\Amat\veceta+\muvec$, then $\veceta$ is a vector with standard a Gaussian distribution. Consequently, $q_{\basisInd}$ can be calculated via the integral of $\veceta$:
\begin{align}
\label{equ:xi2eta}
\nonumber     q_{\basisInd} =& \int_{-\infty}^{\infty} \vecpar^{\basisInd} \mathcal{N}(\vecpar | \muvec, \Sigmamat)   d\vecpar\\
=&\int_{-\infty}^{\infty} (\Amat\veceta+\muvec)^{\basisInd} \frac{\exp(-\veceta^T\veceta)}{\sqrt{(2\pi)^d}} d\veceta.
\end{align}
In this formulation, $(\Amat\veceta +\muvec)^{\basisInd}$ is not the product of univariate functions of each $\eta_i$, therefore, the above integration is still hard to compute. We show that $q_{\basisInd}$ can be computed exactly with an efficient functional tensor-train method.


\subsection{Functional Tensor-Train Implementation}

In this subsection, we show that there exists a matrix $\Gmat_0 \in \mathbb{R}^{1\times r_0}$ and a set of univariate functional matrices $ \Gmat_i (\eta_i) \in \mathbb{R}^{r_{i-1} \times r_i}$ for $i=1,\cdots d$ and with $r_{d}=1$, such that
\begin{equation}
\label{equ:de-correFun}
(\Amat\veceta +\muvec)^{\basisInd} =\Gmat_0 \Gmat_1(\eta_1)\Gmat_2(\eta_2)\ldots \Gmat_d(\eta_d).
\end{equation}
As a result, we have the following cheap computation
\begin{equation}
\nonumber
q_{\basisInd} = \Gmat_0 \mathbb{E}[\Gmat_1(\eta_1)]\mathbb{E}[\Gmat_2(\eta_2)]\ldots \mathbb{E}[\Gmat_d(\eta_d)].
\end{equation}


\subsubsection{Formula  for \eqref{equ:de-correFun}  with $|\basisInd|=1$}
Recall that from $\vecpar=\Amat\veceta+\muvec$, we have
\begin{equation}
\label{equ:xij}
\xi_j=a_{j1}\eta_1+ a_{j2}\eta_2+\ldots+a_{jd}\eta_d+\mu_j, \forall j=1,\ldots,d.
\end{equation}
Here $a_{kj}$ denotes the $(k,j)$-th element of $\mat{A}$.

\textbf{Theorem 3.1}  [Theorem 2, \cite{oseledets2013constructive}]
\textsl{For any function written as the summation of  univariate functions:
$$f(x_0,\ldots,x_d)   =\omega_0(x_0)+\ldots+\omega_d(x_d),$$
it holds that
\begin{align*}
f(x_0,x_1,\ldots,x_d)&
=
\left(\begin{array} {cc}
\omega_0(x_0) & 1
\end{array}\right)
\left(\begin{array} {cc}
 1 & 0\\
 \omega_1(x_1)&1
\end{array}\right)\\
&\ldots
\left(\begin{array} {cc}
 1&0\\
 \omega_{d-1}(x_{d-1})&1
\end{array}\right)
\left(\begin{array} {c}
 1\\
 \omega_d(x_d)
\end{array}\right).
\end{align*}}

Applying Theorem 3.1 to \reff{equ:xij}, we can derive a functional tensor train decomposition for  $\xi_j$:
\small
\begin{equation}
\label{equ:xiTTD}
\xi_j=
\left(\mu_j \ 1\right)
\left(\begin{array} {cc}
 1 & 0\\
 a_{j1} \eta_1&1
\end{array}\right)
\ldots
\left(\begin{array} {cc}
 1&0\\
 a_{j(d-1)} \eta_{d-1}&1
\end{array}\right)
\binom{1}{a_{jd} \eta_d}.
\end{equation}\normalsize
Then the expectation equals to
\begin{equation}
\label{equ:ExiTTD}
\mathbb{E}[\xi_j]=
\left(\mu_j \ 1\right)
\left(\begin{array} {cc}
 1 & 0\\
 0 &1
\end{array}\right)
\ldots
\left(\begin{array} {cc}
 1&0\\
 0&1
\end{array}\right)
\binom{1}{0} =\mu_j.
\end{equation}
The obtained functional tensor trains \reff{equ:xiTTD} can be reused to compute high-order moments.

\begin{algorithm}[t]
\label{alg:moments}
\caption{A functional tensor-train method for computing the basis functions of Gaussian mixtures}
      \SetKwInput{Input}{Input}
      \SetKwInput{Output}{Output}
\Input{The mean value $\muvec_i$, covariance $\Sigmamat_i$ and weight $w_i$ for Gaussian-mixtures, and the  order $p$.}
\For{$i=1,\ldots,n$}
{
Calculate the Cholesky factor $\Amat$ via $\Sigmamat_i=\Amat\Amat^T$;\\
Calculate the functional tensor trains for the first-order and high-order monomials via \reff{equ:xiTTD} and \reff{equ:tt_higherorder}, respectively;\\ Obtain the moments via \reff{equ:ExiTTD}
and \reff{equ:Exi}.
}
Assemble the multivariate moment matrix $\Mmat$ in \reff{equ:mommat};\\
Compute the basis functions via \reff{equ:basis}. \\
\Output{The multivariate basis functions $\{{\multiGPC}_{\basisInd}(\vecpar) \}$.}
\end{algorithm}

\subsubsection{Recurrence Formula  for $1<|\basisInd|\le2p$}

For each $\basisInd$ with $1<|\basisInd|\le2p$, there exist $\basisInd_1$ and $\basisInd_2$ with $|\basisInd_1|,|\basisInd_2|\le p$, such that
\begin{equation*}
\vecpar^{\basisInd} = \vecpar^{\basisInd_1}\cdot \vecpar^{\basisInd_2}, {\rm \ where \ }
\basisInd = \basisInd_1+\basisInd_2.
\end{equation*}
According to \reff{equ:TTkron}, the tensor-train representation of
$\vecpar^{\basisInd}$ can be obtained as the Hadamard
product of the tensor trains of $\vecpar^{\basisInd_1}$ and $\vecpar^{\basisInd_2}$.
Suppose that $\vecpar^{\basisInd_1}=\mathbf{E}_0(\muvec)\mathbf{E}_1(\eta_1)\cdots \mathbf{E}_{d}(\eta_d)$ and $\vecpar^{\basisInd_2}=\mathbf{F}_0(\muvec)\mathbf{F}_1(\eta_1)\cdots \mathbf{F}_{d}(\eta_d)$, then
\small
\begin{align}
\vecpar^{\basisInd}& =
 \Gmat_0(\muvec)\Gmat_1(\eta_1)\cdots \Gmat_{d}(\eta_d),\; \rm{with \ }\Gmat_i = \mathbf{E}_i \otimes \mathbf{F}_i,\ \forall i=0,\cdots,d.
 \label{equ:tt_higherorder}
\end{align} \normalsize
Because $\eta_i$'s are mutually
independent, finally we have
\begin{equation}
\label{equ:Exi}
\mathbb{E}[\vecpar^{\basisInd}] = \Gmat_0(\muvec)\mathbb{E}[\Gmat_1(\eta_1)]\cdots \mathbb{E}[\Gmat_{d}(\eta_d)].
\end{equation}
The moments therefore can be easily computed via small-size matrix  products.

The basis construction framework is summarized in Alg.~\ref{alg:moments}


 \section{AN ADAPTIVE sparse solver}
 \label{sec:sparsesolver}

With the proposed basis functions $\{\multiGPC_{\basisInd}(\vecpar)\}$ for non-Gaussian correlated cases, now we proceed to compute $\{c_{\basisInd}\}$ and express $\out(\vecpar)$ as the form in \reff{eq:ygpc}. The standard stochastic spectral methods~\cite{ghanem1991stochastic,zzhang:tcad2013,xiu2005high} cannot be
directly applied, therefore we improve a sparse regression method via adaptive sampling.

Again, we resort all basis functions and their weights $c_{\basisInd}$ based on the graded lexicographic order, and denote them as $\Psi_j(\vecpar)$ and $c_j$ for $j=1,2,\cdots, N$. Given $M$ pairs of parameter samples and simulation samples $\{ \vecpar_i, \out(\vecpar_i)\}$ for $i=1,2,\cdots, M$, we have a linear equation system
\begin{equation}
\Phimat \cvec = \yvec, \text{\ with\ } \Phi_{ij}=\multiGPC_{j}(\vecpar_i), \ y_i=\out(\vecpar_i),
 \label{equ:linearsystem}
 \end{equation}
where $\Phimat\in\mathbb{R}^{M\times N}$ store the value of $N$ basis functions at $M$ parameter samples. In practice, computing each output sample $\out(\vecpar_i)$ requires calling a computationally expensive device or circuit level simulator. Consequently, it is highly desirable to solve \eqref{equ:linearsystem} when $M\ll N$.

\subsection{Why Do Sparse Solvers Work?}
When $M\ll N$, there are infinitely many solutions to \reff{equ:linearsystem}. A popular method to overcome this issue is to seek for a sparse solution by solving the $\ell _0$-minimization problem
\begin{equation}
\min_{\cvec\in\mathbb{R}^N} \|\cvec\|_0\quad \text{s.t.}\quad \Phimat\cvec=\yvec,
\label{equ:l0optim}
\end{equation}
where $\|\cvec\|_0$ denotes the number of nonzero elements.

The success of an $\ell_0$-minimization relies on some conditions. Firstly, the solution $\mat{c}$ should be really sparse. This is generally true for high-dimensional uncertainty quantification problems. Secondly, the exact recovery of a sparse $\mat{c}$  requires the matrix $\Phimat$ to have the restricted isometry property \cite{candes2006stable}:  there exists a positive value $\delta_s$ such that
\begin{equation}
(1-\delta_s)\|\cvec\|_2^2 \le \|\Phimat \cvec\|_2^2\le(1+\delta_s)\|\cvec\|_2^2
\end{equation}
holds for any $\|\cvec\|_0\le s$. Here, $\|\cdot\|_2$ is the Euclidean norm. Intuitively, this requires that all columns of $\Phi$ are nearly orthogonal to each other. If the $M$ samples are chosen randomly, then for any $k\neq j$, the inner product of the $k$-th and $j$-th columns of the matrix $\Phimat$ is
\begin{equation}
\label{eq:rip_test}
\frac{1}{M}\sum \limits_{i=1}^M \left(\Psi_k(\vecpar_i) \Psi_j(\vecpar_i)\right) \approx \mathbb{E}\left[ \Psi_k(\vecpar) \Psi_j(\vecpar) \right]=0.
\end{equation}
In other words, the restricted isometry property holds with a high probability due to the orthonormal property of our basis functions. Consequently, the formulation \eqref{equ:l0optim} can provide an accurate solution with a high probability.

Once the above conditions hold, \eqref{equ:l0optim} can be solved via various numerical solvers. We employ COSAMP~\cite{needell2009cosamp}, because it can significantly enhance the sparsity of $\mat{c}$.

\begin{algorithm}[t]
\label{alg:AdaptiveSparseSolver}
\caption{An adaptive sparse  solver} 
      \SetKwInput{Input}{Input}
      \SetKwInput{Output}{Output}
\Input{Input a set of candidate samples $\Xi_0$ and  basis functions $\multiGPC_1,\ldots,\multiGPC_N$.}

Choose the initial sample set $\Xi\subset \Xi_0$ via the rank-revealing QR factorization, with $|\Xi|\ll N$.      \\
Call the simulator to calculate $\out(\vecpar)$ for all $\xi \in \Xi$.

\For{Outer iteration $T=1,2,\ldots$}
{Solve the $\ell_0$ minimization  problem \reff{equ:l0optim} to obtain $\cvec$, such that the sparsity $s$ is less than $|\Xi|$.\\
\For{Inner iteration $t=1,2,\ldots,t_{\max}$ }
{
 Fix the indices of the nonzero elements in $\cvec$;\\
Choose a new sample $\vecpar$ from $\Xi_0\setminus \Xi$ by  formula \reff{equ:easy_detmax}, and update the sample set $\Xi=\Xi\cup \{\vecpar\}$;\\
Call the simulator to get $\out(\vecpar)$ at the  sample $\vecpar$;\\
Update the nonzero elements of  $\cvec$ via \reff{equ:updatec}.\\
\If{the stopping criterion is satisfied}{Stop}
}
}
\Output{The basis function coefficient  $\cvec$.}
\end{algorithm}

\subsection{Improvement via Adaptive Sampling}

Previous publications mainly focused on how to solve a under-determined linear system \eqref{equ:linearsystem}, and the system equations were set up by simply using  random samples~\cite{xli2010}. We note that some random samples are informative, yet others are not. Therefore, the performance of a sparse solver can be improved if we select a subset of ``important" samples to set up the linear equation. Here we present an adaptive sampling method: it uses a rank-revealing QR decomposition to pick  ``important" initial samples, and uses a D-optimal criteria to select subsequent samples.

\begin{figure*}[t]
	\centering
		\includegraphics[width=6.6in, height=3.5in]{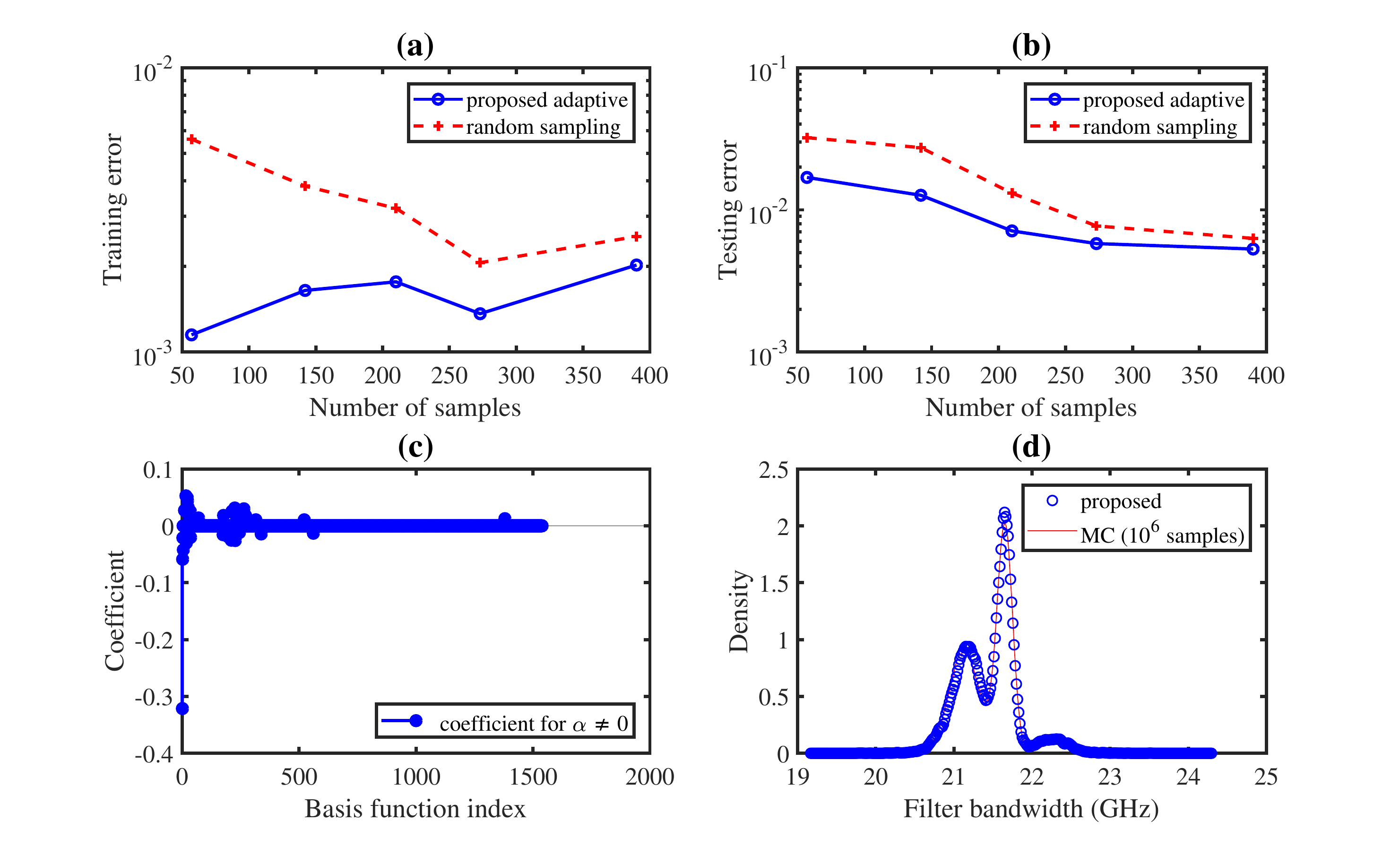}
\caption{Results for the phonic bandpass filter.
(a) training error on 390 samples; (b) testing error on 9000 new samples;
(c) calculated coefficients/weights of our proposed basis functions;
(d) probability density functions of the 3-dB bandwidth obtained with our proposed method and Monte Carlo (MC), respectively.}
	\label{fig:res_PIC19}
\end{figure*}

\textbf{Initial Sampling.} Given a pool of randomly generated candidate samples   for $\vecpar$, we first select a small number of initial samples from this pool. By evaluating all basis functions on the candidate samples, we form a matrix $\Phimat$ whose $j$-th rows stores the values of $N$ basis functions at the $j$-th candidate sample. In order to choose the $r$ most informative rows, we perform a rank-revealing QR factorization~\cite{gu1996efficient}:
\begin{equation}
\Phimat^T\mathbf{P} = \mathbf{Q}\left[\begin{array}{cc}
 \mathbf{R}_{11}&\mathbf{R}_{12}\\
 \mathbf{0}&\mathbf{R}_{22}
 \end{array}\right].
\end{equation}
Here $\mathbf{P}$ is a permutation matrix; $\mathbf{Q}$ is an orthogonal matrix; $\mathbf{R}_{11}\in\mathbb{R}^{r\times r}$ is a upper triangular matrix with nonnegative diagonal elements. Usually,
 $\mathbf{P}$ is chosen such that the minimal singular value of $\mathbf{R}_{11}$ is large, and the maximal singular value of $\mathbf{R}_{22}$ is  sufficiently small.
Consequently, the first $r$ columns of $\mathbf{P}$ indicate the most informative $r$ rows of $\Phimat$. We keep these rows and the associated parameter samples.

\begin{figure}[t]
	\centering
		\includegraphics[width=2.1in]{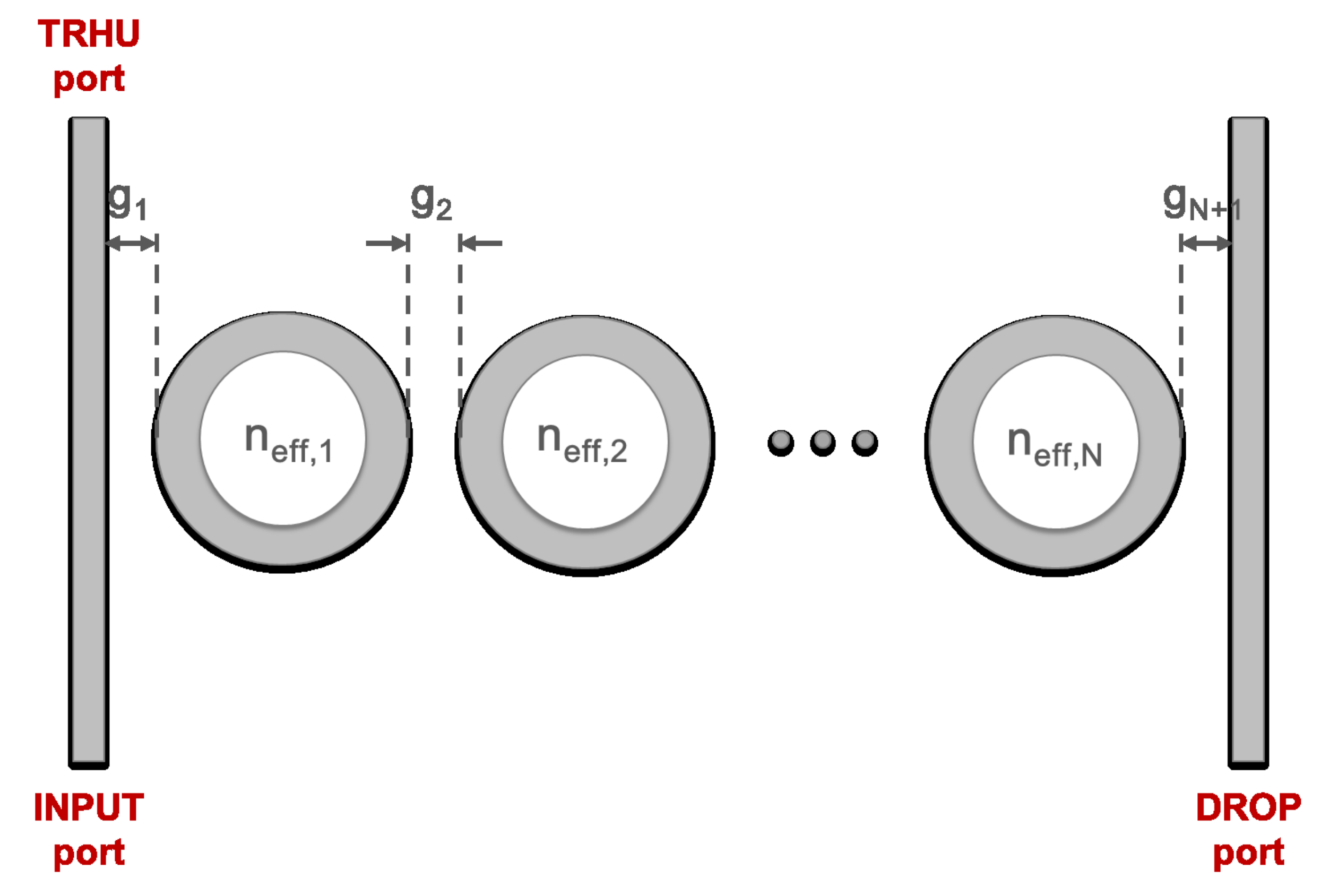}
\caption{A bandpass filter with 9 ring resonators.}
	\label{fig:photonics}
\end{figure}

\textbf{Adding New Samples and Solution Update.} The above initial sampling has generated a $r$-by-$N$ matrix $\Phimat$. Now we further add  new informative samples. Our idea is motivated by the D-optimal sampling update in~\cite{diaz2017sparse}, but differs in the numerical implementation.

Assume that our sparse solver has computed $\mat{c}$ based on the available samples, then we fix the indices of the nonzero elements in $\mat{c}$ and update the samples and solution sequentially.
Specifically, denote the locations of nonzero coefficients  as $\mathcal{S}=\{i_1,\ldots,i_s\}$ with $r>s$, and denote
$\Phimat_1\in\mathbb{R}^{r\times s}$ as the sub-matrix of  $\Phimat$ generated by extracting the columns of $\mathcal{S}$.
 The next most informative sample $\vecpar$ associated with the row vector $\xvec=[\multiGPC_{i_1} (\vecpar) \cdots,\multiGPC_{i_s}(\vecpar)]$ can be found via solving the following problem:
\begin{equation}
\label{equ:detmax}
\max_{\xvec\in \Omega}\quad \det(\Phimat_1^T\Phimat_1 + \xvec^T \xvec),
\end{equation}
where $\Omega$ is value of basis functions in $\mathcal{S}$ for all candidate samples.
In practice, we do not need to compute the above determinant for each sample. Instead, the matrix determinant lemma \cite{harville1997matrix} shows $\det(\Phimat_1^T\Phimat_1+\xvec^T\xvec) = (\det(\Phimat_1^T\Phimat_1))(1+\xvec$ $(\Phimat_1^T\Phimat_1)^{-1}\xvec^T)$, therefore, \reff{equ:detmax} can be solved via \begin{equation}
  \label{equ:easy_detmax}
  \max_{\xvec\in\Omega}\quad \xvec(\Phimat_1^T\Phimat_1)^{-1}\xvec^T.
  \end{equation}

   After getting the new sample,
we update the matrix $\Phimat_1:=\left[\begin{array}{c}\Phimat_1\\ \xvec\end{array}\right]$,  update $(\Phimat_1^T\Phimat_1)^{-1}$ via the Sherman-Morrison formula \cite{sherman1950adjustment}, and recompute the $s$ nonzero elements of $\cvec$ by
\begin{equation}
\cvec_1=(\Phimat_1^T\Phimat_1)^{-1}\Phimat_1^T \yvec.
\label{equ:updatec}
\end{equation}
Inspired by \cite{malioutov2010sequential}, we stop the procedure if $\cvec_1$ is very close to the value of the previous step or the maximal iteration number is reached.
The whole framework is summarized in Alg.~\ref{alg:AdaptiveSparseSolver}.

\section{Numerical Results}
\label{sec:result}

In this section, we validate our algorithms by two real-world examples: a photonic bandpass filter and 7-stage CMOS ring oscillator.
All codes are implemented in MATLAB and run on a desktop with a 3.40-GHz CPU and a 8-GB memory. For each example, we adaptively select a small number of samples from a pool of $1000$ candidate samples,
and we use $9000$ different samples for accuracy validation. Given  random samples, we define the relative error based on \reff{equ:linearsystem}:
\begin{equation}
\label{equ:res}
\epsilon = \|\Phimat \cvec - \yvec\|_2 / \|\yvec\|_2.
\end{equation}
We call $\epsilon$ as a testing error if the samples are those used in our sparse solver, and as a training error if the new set of $9000$ samples are used to verify the predictability.
\begin{table}
\caption{Accuracy comparison on the photonic bandpass filter. The underscores indicate precision.}
\begin{tabular}{c|c|c|c|c}
\hline
method& Proposed &\multicolumn{3}{c}{Monte Carlo} \\ \hline
\# samples & 390 & $10^2$& $10^4$& $10^6$\\ \hline
mean (GHz)& 21.4717 & 2\underline{1}.5297& 21. \underline{4}867& 21.4\underline{7}82\\  \hline
\end{tabular}
\label{tab:pic}
\end{table}

\subsection{Photonic Bandpass Filter (19 Parameters)}

\begin{figure*}[t]
	\centering
		\includegraphics[width=6.6in, height=3.8in]{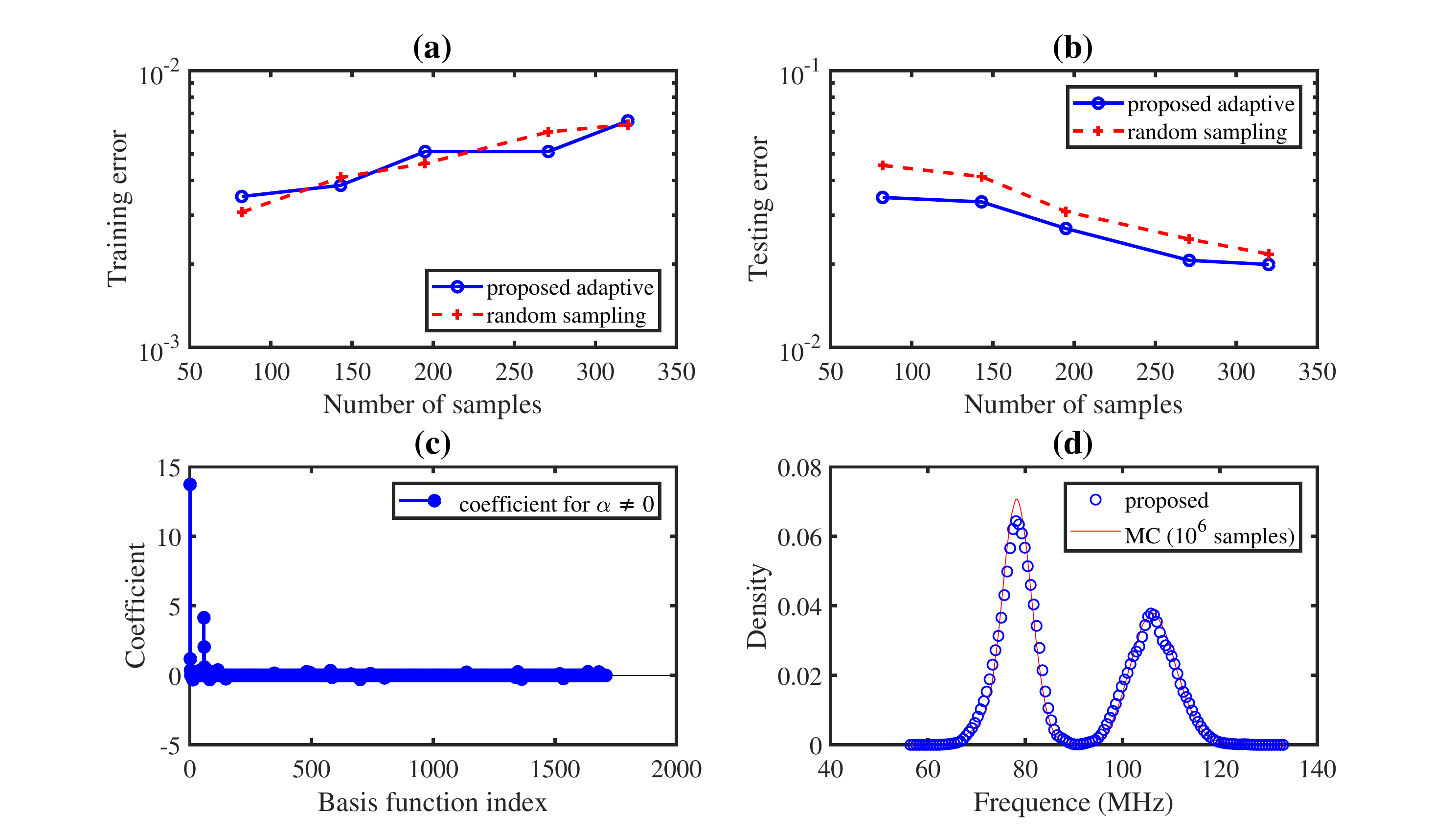}
\caption{Results for the CMOS ring oscillator. (a) training error on 320 samples; (b) testing error on 9000 new samples;
(c) calculated coefficients/weights of our proposed basis functions;
(d) probability density functions of the oscillator frequency obtained with our proposed method and Monte Carlo (MC), respectively.}
	\label{fig:res_ring57}
\end{figure*}

Firstly we consider the photonic bandpass filter in Fig.~\ref{fig:photonics}. This photonic IC has 9 ring resonators, and it was originally designed to have a 3-dB bandwidth of 20 GHz, a 400-GHz free spectral range, and a 1.55-$\mu$m operation wavelength. A total of $19$ random parameters are used to describe the variations of the effective phase index ($n_{\rm eff}$) of each ring, as well as the gap ($g$) between adjoint rings and between the first/last ring and the bus waveguides. These non-Gaussian correlated random parameters are described by a Gaussian-mixture joint probability density function. 
\begin{figure}[t]
	\centering
		\includegraphics[width=2.1in]{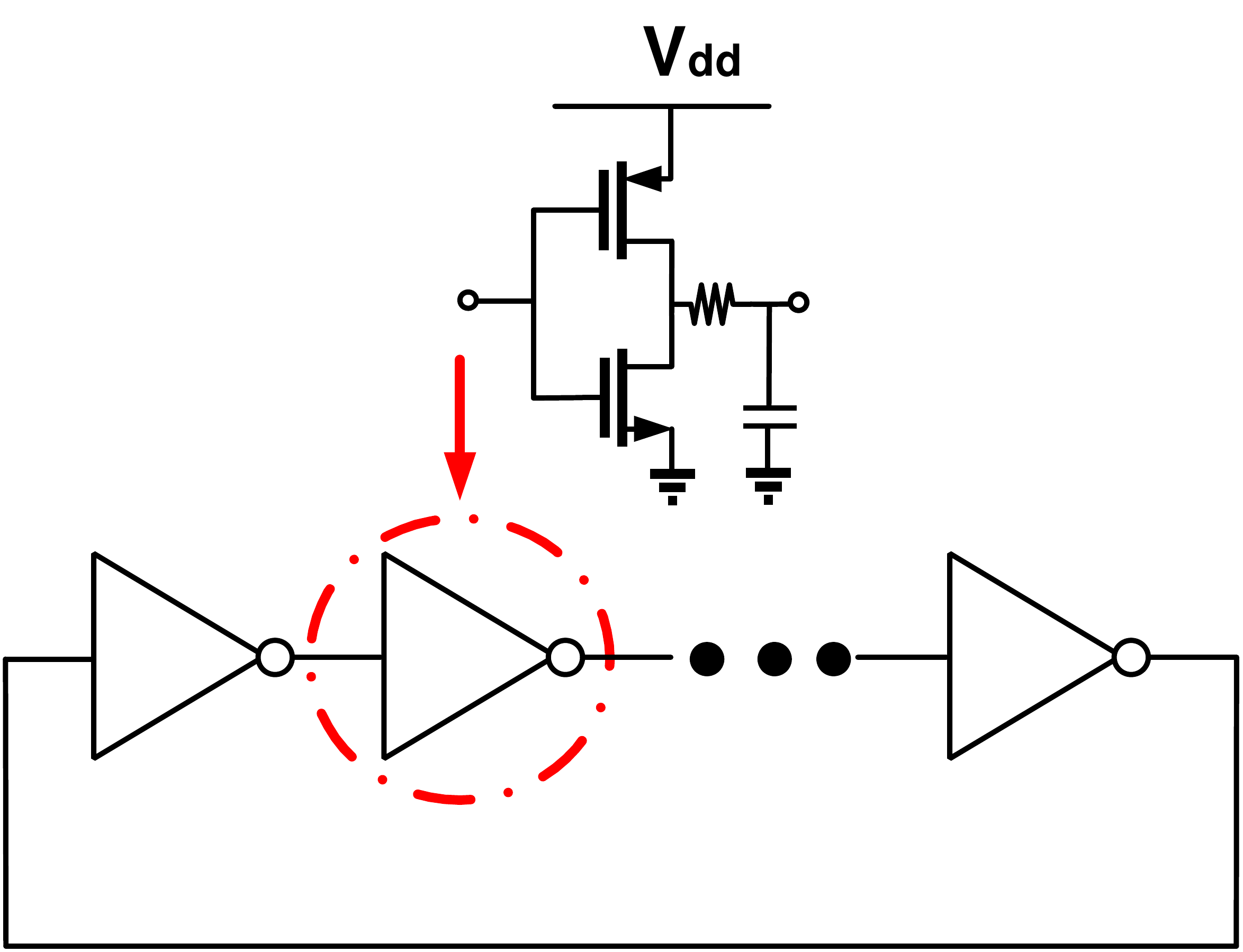}
\caption{Schematic of a CMOS ring oscillator. }
	\label{fig:ring}
\end{figure}

We approximate the 3-dB bandwidth $f_{3 {\rm dB}}$ at the DROP port using our proposed basis functions with a total order bounded by $p=3$. The numerical results are shown in Fig.~\ref{fig:res_PIC19}. Fig.~\ref{fig:res_PIC19} (b) clearly shows that
the adaptive sampling method leads to significantly lower testing (i.e., prediction) errors when a few samples are used, because it chooses important samples. Finally, we use 390 samples to assemble a linear system and solve it by $\ell_0$ minimization, and obtain the sparse coefficients of our basis functions in Fig.~\ref{fig:res_PIC19} (c). Although a third-order expansion involves more than $1000$ basis functions, only a few dozens are important.
Fig.~\ref{fig:res_PIC19} (d) shows the predicted probability density function of the filter's 3-dB bandwidth, and it matches well with the result from Monte Carlo. More importantly, it is clear that our algorithm can capture accurately the multiple peaks in the output density function, and these peaks can be hardly predicted using existing stochastic spectral methods.

In order to demonstrate the effectiveness of our method in detail, we compare the computed mean value of $f_{3 {\rm dB}}$ from our methods and from Monte Carlo in Table~\ref{tab:pic}. Our method provides a closed-form expression for the mean value. Monte Carlo method converges very slowly, and requires $2564\times$ more simulation samples to achieve the similar level of accuracy (with 2 accurate fractional digits).

\subsection{CMOS Ring Oscillator (57 Parameters)}

We continue to consider the 7-stage CMOS ring oscillator in Fig.~\ref{fig:ring}. This circuit has $57$ random parameters describing the variations of threshold voltages, gate-oxide thickness, and effective gate length/width. We use Gaussian mixtures to describe the strong non-Gaussian correlations.


We use a 2nd-order expansion of our basis functions to model the oscillator frequency. The simulation samples are obtained by calling a periodic steady-state simulator repeatedly. The detailed results are shown in Fig.~\ref{fig:res_ring57}. Similar to the previous example, our adaptive sparse solver produces a sparse and highly accurate stochastic solution with better prediction behaviors than the standard compressive sensing does. The proposed basis functions can well capture the multiple peaks of the output probability density function caused by the strong non-Gaussian correlation.

Table~\ref{tab:ring} compares our method with Monte Carlo. Our method is about $3125\times$ faster than Monte Carlo to achieve a precision of one fractional digit   for the mean value.

\begin{table}
\caption{Accuracy comparison on the CMOS ring oscillator. The underscores indicate precision.}
\begin{tabular}{c|c|c|c|c}
\hline
method& Proposed &\multicolumn{3}{c}{Monte Carlo} \\ \hline
\# samples & 320 & $10^2$& $10^4$& $10^6$\\ \hline
mean (MHz) &90.5441 & 89.7795& 9\underline{0}.4945& 90.\underline{5}253\\
\hline
\end{tabular}
\label{tab:ring}
\end{table}
\section{Conclusion}
\label{sec:conclusion}
This paper has presented some theories and algorithms for the fast uncertainty quantification of electronic and photonic ICs with non-Gaussian correlated process variations.
We have proposed a set of  basis functions for non-Gaussian correlated cases. We have also presented a functional tensor train method to efficiently compute the high-dimensional basis functions.
In order to reduce the computational time of analyzing process variations,
we have proposed an adaptive sampling sparse solver, and this algorithm only uses a small number of important simulation samples to predict the uncertain output. The proposed approach has been verified with two electronic and photonic ICs. On these benchmarks, our method has achieved very high accuracy in predicting the multi-peak output probability density function and output mean value. Our method has achieved $2500\times$ to $3100\times$ speedup over Monte Caro to achieve the similar level of accuracy. To our best knowledge, this is the first non-Monte-Carlo uncertainty quantification approach that can handle non-Gaussian correlated process variations without any error-prone de-correlation steps. Many novel algorithms can be further developed based on our results.


\section*{Acknowledgment}
This work was supported by NSF-CCF Award No. 1763699, the UCSB start-up grant and a Samsung Gift Funding.

\balance

\bibliographystyle{IEEEtran}
\small{
\bibliography{RefList}
}

\end{document}